\documentclass[12pt]{amsart}
\usepackage{geometry}                
\usepackage{graphicx}
\usepackage{amssymb}
\usepackage{epstopdf}
\newtheorem{definition}{Definition}
\newtheorem{axiom}{Axiom}

\title{Complete Totalities}
\author{Rafi Shalom}
\email{rafi.shalom@gmail.com}

\begin{document}

\bibliographystyle{plain}

\maketitle


\begin{abstract}

The cumulative hierarchy conception of set, which is based on the conception that sets are inductively generated from ``former'' sets, is generally considered a good way to create a set conception that seems safe from contradictions. This imposes two restrictions on sets. One is a ``limitation of size,'' and the other is the rejection of non-well-founded sets. Quine's NF system of axioms, does not have any of the two restrictions, but it has a formal restriction on allowed formulas in its comprehension axiom schema, which reflects a similar notion of elements being prior to sets. Here we suggest that a possible reason for set antinomies is the tension between our perception of sets as entities formed from elements by an imaginary aggregation operator, and our wish to regard sets as existing ``at once.'' A new approach to sets as totalities is presented based on a notion of ``concurrent aggregation,'' which instead of avoiding ``viscous circles,'' acknowledges the inherent circularities of some predicates, and provides a way to characterize and investigate these circularities.


\end{abstract}


\section{Introduction}

Naive set theory is known to be contradictory. The reaction to set theoretic paradoxes appears to have been restrictive in the sense that it has been assumed that the contradictions are a result of assuming too much. This is a natural response. However, by examining the restrictions imposed by two different foundational approaches and their philosophical justifications, it would be possible to offer a new interpretation of set theoretic paradoxes. This interpretation asserts that the paradoxes are a result of the tension between two equally fundamental aspects of sets. The first aspect is related to the way we perceive sets and their elements, and the need to ``mentally reach'' sets from their elements, and vice versa. This suggests an implicit notion of precedence on sets which conflicts with another aspect we may wish to keep in our mathematics - having a universe of sets free from hidden or explicit notions of precedence.

This observation promotes the notion of concurrent aggregation, and leads to a new conception of set, which is similar to the naive conception, yet one that is an expansion of the conception of set rather than a restriction of it. In this new conception, instead of trying to exclude circularities within the properties that define their respective totalities\footnote{By ``totality'' we mean the collection of all objects that satisfy some predicate.}, these circularities are allowed to surface. The philosophical aspects of this proposition allow several insights about the way we perceive sets, and the mathematical aspects provide an interesting new universe of sets which includes a variety of sets required for ``everyday mathematics,'' and might allow a different perspective on several set theoretic issues.

In section \ref{sec:cumu} we examine two restrictions that occur in cumulative hierarchies, namely ``limitation of size'' and the rejection of non-well-founded sets (NWFS will be used here to abbreviate ``non-well-founded sets''), which are sometimes also referred to as ``hypersets.'' Non-well-founded sets are sets with infinitely descending $\in$-sequences, such as sets that contain themselves as elements. Cumulative hierarchies may arise from the notion of an inductively generated set, or from the notion that elements of a set are ``prior to'' it, in the sense of ontological dependence. Both motivations are discussed in Boolos' article on the iterative conception of set \cite{BoolosIterativeConceptionOfSet}, which is sometimes seen as a justification for ZF. 

Even though the rejection of NWFS appears to be a less arbitrary restriction than ``limitation of size,'' we try to gain a better perspective of it by examining it through three different angles. The first is the existence of well-established theories of NWFS, and their justifications. A second angle is the philosophical justifications provided by Boolos for cumulative hierarchies. These justifications appear to stem from an interpretation of sets based on the notion that elements of a set are ``prior to'' it. The third perspective tracks down the conceptual differences between stage or rank theories, such as Boolos' iterative conception, and ZF which is assumed to be motivated by these conceptions. While the assertion that all sets are well-founded is a logical consequence of rank theories, and an integral part of them, in ZF it is possible to exclude the regularity axiom (sometimes also called the foundation axiom) which states that all sets are well-founded. We provide a simple motivation for the axioms of ZF which stems from an immediate interpretation of the axioms, that even in its strongest form does not imply the foundation axiom. We conclude that ZF, unlike the iterative conception, is not necessarily based on the notion that elements of sets are ``prior to'' them. 

In section \ref{sec:nf} Quine's NF is discussed. NF does not include neither limitation of cumulative hierarchies. Indeed Quine expresses concern about Russell's interpretation of the paradoxes, which limits our ability to define sets by quantifying over all sets. However, Quine's argumentation for NF is not alien to Russel's type theory. It is also explicitly mindful of the paradoxes, and the formal/syntactic limitation Quine imposes on formulas in his specification axioms, namely ``stratification,'' is based on similar notions of precedence. One of the possible motivations for NF, which enhances the iterative conception with the possibility of aggregating sets from later stages, and introduces aggregations that uses ``placeholders'' which results in circular membership, shows that NF also has much in common with the iterative conception. 

Finally in section \ref{sec:complete} we present the notion of a complete totality, which applies a form of aggregation we call ``concurrent aggregation'' to sets identified by a predicate. Concurrent aggregation uses a form of aggregation with placeholders, but it is free of stages, stratification, and other layered notions. The theory is characterized in an informal way, and immediate consequences are reached.

\section{Cumulative Hierarchies}
\label{sec:cumu}

One reaction to the paradoxes, fairly immanent and understandable, leads to a view of sets as being ``piled'' in layers, each layer is ``built from'' sets from former layers, which seems to demand two kinds of limitations on sets. The first limitation is sometimes called ``limitation of size,'' and the other is the exclusion of NWFS. There are many specific formulations and axiomatizations of theories that assert these limitations (Russel's theory of types, axioms systems of Morse-Kelley, and of Von Neumann to name a few of the least recent ones). These limitations occur in the widely used ZF and ZFC systems of axioms. Several philosophical justifications have been offered to support such theories. A defensive motivation for a layered view of mathematics is expressed by Russell \cite{Rus:08}. Russell blames self-referential aspects of the paradoxes as their cause. Boolos provides arguments for the rejection of NWFS, and also positive incentives for a hierarchical universe of sets as part of his iterative conception of set \cite{BoolosIterativeConceptionOfSet, boolos:1989}. This conception is offered as philosophical support for ZF, which for Boolos and others is synonymous with set-theory. Since Boolos seems to have the most systematic philosophical justification of this view of sets, we will be mostly concerned with his argumentation in this section.

\subsection{Limitation of Size}

Cumulative hierarchies reject the notion of a universal set, that is, a set that contains all sets. A universal set is not in itself contradictory, and there are theories of sets that contain a universal set, such as Qunie's NF \cite{QuineNewFoundationsForMathematicalLogic}. However, a universal set cannot exist if we expect that every set has a power set, and that the power set of each set has greater cardinality. In all versions of cumulative hierarchies we always keep the latter assertion, and usually keep the former assertion. Worse yet, in a cumulative hierarchy approach, a universal set usually leads directly to the existence of Russel's set by an application of a comprehension axiom. So in a cumulative hierarchy there is a need to cut the hierarchy at some point and declare that from that point on, the objects are ``classes'' rather than sets, or ``too big'' to be sets. 

As a doctrine, ``limitation of size'' determines which objects are sets, and which are not, according to their size relative to that of other objects in a specific universe of sets. In \cite{boolos:1989} Boolos provides, among other things, a formulation of ``limitation of size'' in two different ways, a stronger one and a weaker one. The only difference is that in the weaker version we might not want to allow as sets some objects that do not have a bijection with the class of all sets. He also provides a version of Frege's system with a modified rule V which avoids a contradiction by allowing only ``small concepts.'' 

This analysis of ``limitation of size'' by Boolos, allows him to deduce several aspects of ZF that could not be accounted for by his iterative conception of set \cite{BoolosIterativeConceptionOfSet}. For example, Boolos feels that the iterative conception can support stages (the equivalent of ranks) up to the first uncountable ordinal or even only up to the first nonrecursive ordinal, so he originally has to support the replacement axioms in ZF with a practical consideration: ``... the reason for adopting the axioms of replacement is quite simple: they have many desirable consequences and (apparently) no undesirable ones.'' Later \cite{boolos:1989} he fills up this philosophical gap by noting that ``limitation of size'' can provide what the iterative conception cannot, ending with the conclusion that: ``Perhaps one may conclude that there are at least two thoughts `behind' set theory.''

However interesting ``limitation of size'' can be as a doctrine, it is indifferent to the actual scope of the cumulative hierarchy, which leaves the question of ``how big is too big'' open. There are many possible considerations, some more philosophically grounded, and some are practical. Even some of the writers who aim to provide a motivation for set theory, may not see the importance of argumentation about its preferred scope. For example, Scott \cite{ScottAxiomatizingSetTheory} provides an axiomatization of set theory that aims to support ZF, and adopts a reflection axiom schema, which is fairly strong, citing only practical considerations.

It would be impossible to survey all existing scope proposals and their justifications, but even without it, it is hard to imagine a good defense of any such decision. For example, even if we are convinced that it is philosophically sound to use one of the highly restrictive schemes, so that we are safer and still able to carry out most of our mathematics, how can this be a foundations of mathematics when it excludes perfectly good mathematics concerning large cardinals? The various restrictions that end with this or that large cardinal seem also to be somewhat arbitrary. It looks as though the only very ``thick line'' one can draw is between the finite and the infinite, but very few would settle for a mathematics without the infinite. Thus while trying to avoid set theoretic antinomies, our defense leads us to a scope of sets which is either unjustified, too restrictive, or a combination of both. 

This apparent foundational dead end may promote skepticism towards a foundational approach to mathematics, or make one want to reject sets as a foundational path. The set conception to be described in section \ref{sec:complete} does not demand scope decisions, and it supports a universal set.

\subsection{Non-Well-Founded Sets}

Boolos is aware that ``Unlike the naive and the iterative conception, {\em{limitation of size}} (in either version) is not a natural view, for one would come to entertain it only after one's preconceptions have been sophisticated by knowledge of set theoretic antinomies ...'' However it remains to be seen to which extent can Boolos' stage theory be called natural, and if another form of restriction it imposes, namely the rejection of NWFS, is as justifiable and free from being a reaction to the paradoxes as Boolos thinks. 

Unlike the ``how big is too big'' issue, asserting that only well-founded sets can be considered sets is generally accepted, and considered natural. The situation is less vague because there is a distinct line that separates well-founded sets from non-well-founded sets. Stage or rank theories such as Boolos', with the set conception they provide, are usually taken as good philosophical justification for the rejection of NWFS, and as a basis for ZF.  The ZF system of axioms contains the axiom of regularity, which rejects the existence of NWFS. Thus sets are assumed to be well-founded in any usual context. 

On the other hand, it is relatively well known that a theory of NWFS is possible and even useful. Peter Aczel \cite{Aczel} provided a unified framework for dealing with NWFS, based on slight modifications of ZF. The regularity axiom is the only one that has to be rejected, and several anti-foundation axioms are offered that differ only by the strength of the criteria for set equality.

Aczel devised a way to describe sets as ``decorations'' of directed graphs in which one of the vertices represents the set, and edges describe membership. Vise versa, graphs can be a ``picture'' of sets, well-founded or otherwise. This graph description of sets makes it easy to describe NWFS, sometimes in a finite way. This is though NWFS appear to be infinite objects because they always have an infinite $\in$-sequence. For example, The non-well-founded set which is a singleton containing itself, which Aczel denotes by $\Omega$, can be described as the decoration of a graph composed of a single vertex with a self loop. There are many other graphs that are pictures of this set, some of them are infinite, but the fact that NWFS should not necessarily be seen as infinite objects did not escape Aczel: ``But a moment's thought should convince the reader that $\Omega$ is a finite an object as one could wish. After all it does have a finite picture.''

Equality for NWFS is a more delicate issue than for well-founded sets, because sets can be members of themselves (or members of members etc.) which obviously influences the usual extensionality rule for set equality. The usual extensionality rule is not only unable to equate two ``different'' singletons containing themselves. Some relatively complicated looking sets such as the ones depicted by the equations $\Omega_1 = \{\Omega_2, \Omega_3\}, \Omega_2 = \{\{\Omega_2\}\}, \Omega_3 = \{\Omega_2\}$ may also be equated to $\Omega$. A unified framework is presented by Aczel to handle his own suggestions and previous work on the subject of NWFS. The framework allows a systematic description of several anti-foundation axioms which differ by the strength of the identity criteria they impose. This aspect of NWFS will be revisited in the set conception we present in section \ref{sec:complete}, for which such set identity considerations are important.

The challenge posed by a systematic and useful theory of NWFS to the cumulative hierarchy notion of a set is stated clearly in a foreword to Aczel's book written by Jon Barwise: ``Aczel's work was motivated by work of Robin Milner in computer science modeling concurrent processes. The fact that these processes are inherently circular makes them awkward to model in traditional set theory, since most straightforward ideas run afoul of the axiom of foundation.'' But the resulting theory by Aczel which adopts NWFS while using standard techniques from set theory results in a linguistic problem: ``However, there is a serious linguistic obstacle to this work, arising out of the dominance of the cumulative conception of set. Just as there used to be complaints about referring to complex numbers as numbers, so there are objections to referring to non-well-founded sets as sets. While there is clear historical justification for this usage, the objection persists and distracts from interest and importance of the subject.''

Barwise's comparison of NWFS with the extension of the conception of number to include complex numbers does less justice to NWFS than they deserve. NWFS should not be seen as legitimate because they are a way to widen the conception of set, which is a usual mathematical practice well described by Meir Buzaglo \cite{Buzaglo2002-BUZTLO}. NWFS were rather rejected as a result of consistency concerns raised by set antinomies. Mathematical objects such as the universal set, and the naive conception of set, are earlier than the cumulative conception. This is probably what Barwise means by a ``clear historical justification'' for NWFS. But even if NWFS were ``just'' a new mathematically well established expansion with useful applications, we would still be inclined to let NWFS be seen as perfectly good sets, instead of having to awkwardly skip between the terms ``set'' and ``hyperset'' whenever hypersets are deemed useful. 

\subsection{Contemplating the Case for Sets Being Strictly Well-Founded}
There are several good reasons why the cumulative conception of set has become so dominant. Even without looking up the justifications for it in the literature, one can clearly appreciate the advantages of a theory of sets that arises from a series of creations of sets from former sets, and it is natural that one would adopt such a view given the knowledge of set antinomies. The antinomies seem to warn us from circularity in our membership relation, and ``constructing'' sets layer by layer from former ones accommodates a reassuring feeling of consistency, because we do not expect that the acceptance of the ``new'' sets we ``produce'' can create a contradiction. The inductive presentation is also similar to the way we usually think about the natural numbers. We concentrate on Boolos' argumentation on the subject that seem to reflect these sentiments.

Concerning argumentation for well-foundedness of sets there can be two basic attitudes - a rejection of NWFS, and claiming that a set conception that includes only well-founded-sets is preferable. Both of them are expressed by Boolos. Boolos argues against NWFS in \cite{BoolosIterativeConceptionOfSet} by saying that ``It is important to realize how odd the idea of something's containing itself is.'' Admittedly, the idea of an $x$ such that $x\in x$ is a bit odd, but can this make a good philosophical argument? After all, all of mathematics is a bit odd, including usual set theoretic concepts Boolos accepts such as the empty set, singletons, and the axiom of choice. Boolos reminds us in a different context \cite{boolos:1989} how odd the idea of a set with a single element, and it being unequal to the element within it, may seem to a non-mathematician. Some of the oddest ideas in mathematics are also some of the most powerful ones. Are infinitesimals and $\sqrt {-1}$ less odd than a universal set? Ultimately mathematicians are concerned with the usefulness of a theory, not its likelihood. The problematic concepts can always be removed. Just like the way the complex numbers can be described as a field of ordered pairs of real numbers, thus excluding $\sqrt {-1}$ from the theory, NWFS can be represented by graphs, thus eliminating the concept of circular membership.

Boolos is aware that ``There does not seem to be any argument that is guaranteed to persuade someone who really does not see ... that these states of affairs are peculiar." Much less this would affect someone who agrees that this may be peculiar, yet does not see a major problem.

The other side of the same coin is Boolos' assertion that the iterative conception he describes later using a ``rough description'' of sets being accumulated in stages, and with its formalization with stage theory, is ``natural.'' This is supposed to validate the theory as one that has its own justification regardless of the paradoxes. Boolos' argumentation is as follows: since naive set theory is contradictory ``one might come to believe that any decision to adopt a system of axioms about sets would be {\em{arbitrary}} ... if consistent, its consistency would be due to certain provisions that were laid down for the express purpose of avoiding the paradoxes that show naive set theory inconsistent, but that lack any independent motivation.'' Contrary to that, Boolos claims this the iterative conception ``often strikes people as entirely natural, free from artificiality, not at all ad-hoc, and one they might perhaps have formulated themselves.'' 

Boolos may be right about the conception being natural, but what could be the implications of having a natural conception of a mathematical entity? Does it mean that this is the only way we can interpret it, or rather that investigating it according to that conception has merit? There is no dispute that limiting our attention to well-founded sets, and working within a framework that ensures that all sets have a rank is convenient and beneficial in many cases. This means that choosing to work in a universe of well-founded sets is legitimate, but it does not mean that our set conception has to exclude NWFS. The natural numbers are probably the most natural, ``not at all ad-hoc'' and ``easy to formulate'' mathematical concept one can imagine, and their theoretical value as a subclass of numbers cannot be denied, yet the real numbers are as much numbers as the natural numbers.

It could be the case that Boolos did not intend to claim that the iterative conception should be regarded as the only way to conceptualize set theory, but only as a way to motivate ZF. However, he does explain his view of how the term ``set'' should be understood. In the beginning of that article Boolos discusses Cantors' definitions of the term. Boolos thinks Cantors' definitions are vague, but for him (justifiably or not is a different question) the definitions ``... do suggest - although, it must be conceded, only very faintly - two important characteristics of sets: that a set is `determined' by its elements ... , and that, ... the clarification of which is one of the principle objects of the theory of the rationale we shall give, the elements of a set are `prior' to it." No wonder later Boolos explains the alleged strangeness of NWFS by saying that ``For when one is told that a set is a collection into a whole of definite elements of our thought, one thinks: Here are some things. Now we bind them into a whole. {\em{Now}} we have a set.'' 

As we shall see in the next section, the exclusion of NWFS is a logical consequence of this notion of priority induced by the membership relation, though it seems Boolos was unaware of that at the time. There are also reasons to believe that the iterative conception and ZF have completely different motivations, which makes the suggestion that the iterative conception motivates at least a part of ZF questionable.

 \subsection{Regularity in ZF and in Stage Theory}
 
If we reject Boolos' philosophical arguments that lend support to regularity, it may still be possible to support regularity as an integral part of ZF by noticing that regularity logically follows from rank theories such as the iterative conception, and adding the assumption that the iterative conception is a justification of the parts of ZF that include regularity. 

Regularity is indeed an integral part of the iterative conception and other rank theories. The fact that regularity is implied by stage theory initially escaped Boolos. This made him introduce set theoretic induction\footnote{Set theoretic induction means that for every property $P$ the statement $\forall x ((\forall y y\in x\rightarrow P(y)) \rightarrow P(x)) \rightarrow \forall x P(x)$ holds. For the version in Boolos' stage theory we first say that a stage is covered by P if P holds for all sets formed at that stage. Then set theoretic induction is: if each stage is covered by P provided all lower stages are covered by P, then all stages are covered by P.} in his original presentation of stage theory. The argument for the adoption of the induction axiom was that the conception indicates that sets are inductively generated. After laying out the axioms of stage theory Boolos adds that ``There is still one important feature contained in our rough description that has not yet been expressed in the stage theory: the analogy between the way sets are {\em{inductively generated}} by the procedure described in the rough statement and the way the natural numbers 0,1,2, ... are inductively generated ...'' Thus for Boolos,  regularity must follow from the iterative conception.

In theories that make stages or ranks of sets explicit, and have a comprehension axiom schema that quantifies over sets of limited rank, it is possible to deduce set-theoretic induction. Such deductions appear in a more compact version of Boolos' stage theory \cite{boolos:1989}, in a theory by Dana Scott that uses membership of ``partial universes'' for ranking \cite{ScottAxiomatizingSetTheory}, and in an axiomatization by James Van Aken \cite{Aken86}, which includes an axiom that defines a two place relation that supports the notion of ``presuppositional strength'' of sets, which is similar to ranks. Set theoretic induction is an essential part, and not merely a logical consequence of such systems. It is easy to notice this in Van Aken's system that has two basic axioms (the above mentioned axiom that defines a relation for set comparison according to their rank, and a comprehension axiom schema that uses this relation) responsible for the derivation of all the usual ZF axioms Boolos' stage theory supports. Remove one of the two axioms and the system collapses. This ``all or noting'' situation occurs also in the other two axiomatizations.

If we accept Boolos' thesis that stage theory supports ZF without replacement and extensionality (and without choice if ZFC is to be justified), accepting regularity as an essential part of ZF is unavoidable. Actually this would mean that we also have to accept the power set axiom, which is derivable from such axiom systems, provided there is no maximal rank\footnote{Actually it is possible to create a Van Aken style system in which the two place relation is relaxed not to indicate strict rank, and the power set axiom follows directly from such a comprehension axiom, but this is beyond the scope of this article.}. To show that this is not the case there is a need to explain what makes regularity ``stick out'' in ZF and ZFC, and to trace the ``thought'' behind ZF which is ``simple'', ``natural'', ``not at all ad-hoc'' from which the axioms supported by stage theory can be deuced, except regularity.


A careful look at ZF and ZFC should help us find a common theme. It is noticeable that most of the ZF axioms are assertions about set existence, usually given the existence of other sets. Most of the axioms can be interpreted to say something like: ``given this or that set (or sets) there is a set such that ...'' or ``there is an infinite set.'' There are only two exceptions. One is extensionality which equates sets that have the same elements, and has a special status in various axiomatizations of set theory. The other exception is regularity which imposes a condition on sets, and is actually the only axiom in ZF that disqualifies sets. Thus it would be fair to say that in ZF (and also in ZFC since the axiom of choice also adheres to the above pattern of axioms of ZF), regularity is rather the ``odd man out.''

How can we characterize the ``thought'' behind this pattern? Unlike the situation in the ``rough description'' that Boolos uses in order to motivate the iterative conception, these sets are not simply built using ``previous'' sets, so the movement is symmetric - not only from elements of a set towards the set, but also from sets towards their elements. This means that given a set we ask ourselves which sets can be ``reached'' or ``constructed'' given other sets. The following simple notion of set construction follows: for all $x$, $x$ is 0-constructible from $x$. $y$ is $n$-constructible from $x$ if every $z\in y$ is up to $(n-1)$-constructible from $x$, or if there is a $z$ which is up to $(n-1)$-constructible from $x$ and $y\in z$ holds. In other words, if we can we get from $x$ to $y$ within $n$ aggregations and disaggregations of sets we already have, $y$ is $n$-constructible from $x$. We are allowed to aggregate and/or disaggregate all the sets we have at once at each stage. For example, getting the union set of a certain set requires (at most) disaggregating the set to its elements, doing the same for the elements, and aggregating what we got. A union set of $x$ is thus always 3-constructible from $x$. Note that according to the above definition, the empty set is 1-constructible from any set, because all the elements of the empty set are 0-constructible from any set, which might not coincide with the original intuition, but this is not important for the point we are about to make.

Even without formalizing this concept it is trivial that the regularity axiom cannot be deduced from a comprehension axiom scheme that quantifies over sets that are $n$-constructible from other sets. Even if we take a very strong version of this axiom that allows starting from more than one set, and quantifying over anything that is finitely constructible\footnote{A formalization of ``finitely constructible from $x$'' does not require a notion of a natural number, and it is easy to formulate in the language of set theory using first order logic with a single type. Since the technical details are fairly straightforward, we refrain from such a detailed exploration here.} (that is, $n$-constructible for some natural $n$) from them, we can prove many ZF axioms - unordered pairs, power set, union, comprehension, and even infinity - but not regularity. This would be easily provable by deducing the version of comprehension that quantifies over finitely constructible sets from ZF without regularity, plus the known fact that ZF without regularity does not imply regularity. Indeed, this form of creating sets from other sets does not seem to be at odds with the existence of, say, the set $\Omega$ mentioned earlier. Actually, this concept of set creation does not seem to be at odds with one of Cantors' definitions for sets that Boolos quoted: ``any collection ... into a whole of definite, well-distinguished objects ... of our intuition or thought.'' A definition is expected to be a two way street. If we know the set when we know its elements, shouldn't we know the elements when we know the set?

If this analysis of the conceptual framework of ZF is true, there should be more suspicion towards assumption that regularity is a part of ``a thought behind ZF''. The assumption that the cumulative approach is essential to ZF is taken for granted by Scott, for example, when he lays out his axiomatization of set theory \cite{ScottAxiomatizingSetTheory}: ``Now Russell made his types {\em{explicit}} in his notation, and Zermelo left them {\em{implicit}}. It is a mistake to leave something so important invisible, because so many people will misunderstand you.'' Perhaps the misunderstanding belongs to Scott, and Zermelo never considered hierarchical notions as a key feature of his set theory. 

This does not mean that rank theories are not a good way to describe a natural conception of set. Boolos' rough description is similar to a description of the natural numbers using zero and the successor operator. Boolos' stage theory relates to his advocated set conception the way Peano's axiom system relates to the natural numbers. The result is a much required description of a subclass of sets. 

To summarize, rank theories essentially include the idea that sets are inductively generated, and make ranks of sets explicit in their axiomatic development. As a result, regularity and other set theoretic consequences are derivable. There is no way to exclude regularity in such axioms systems without rendering them meaningless and dysfunctional. Conversely, in ZF regularity is detachable. This is probably because ZF is an axiom system motivated by the idea that the role of set theoretic axioms is to describe the universe of sets by separate assertions about the existence of sets, thus unlike what Boolos thinks, rank theories do not reveal ``a thought behind ZF,'' but they are a completely different way to perceive sets, which is much more restrictive in some respects.

\section{Quine's New Foundations}
\label{sec:nf}

There is one long-standing and relatively well known axiomatization of sets that does not contain either of the two restrictions mentioned in the previous section. Quine's New Foundations (NF) \cite{QuineNewFoundationsForMathematicalLogic} axiom system is simply an extensionality axiom with a comprehension axiom schema that allows the totalities of predicates defined by formulas in the language of set theory which are ``stratified.'' Stratification demands that the variables on both sides of the $\in$ relation in subformuals of the formula can be enumerated in a non circular way. Specifically, there has to be a function $f$ from the variables of the formula to the natural numbers such that if $x\in y$ is a subformula $f(x) + 1 =  f(y)$ holds, and if $x = y$ is a subformula $f(x) =  f(y)$ holds. Stratified formulas are similar to formulas in Russell's theory of types, without explicit types.

Since the property $x=x$ is allowed, a universal set exists, and since $x\notin x$ is not allowed, Russell's paradox is avoided. The universal set in NF is both an example of allowing ``big'' sets, and of allowing non-well-founded sets. This means that the ``ontology'' derived from NF and from hierarchical theories is very different. Unsurprisingly each side argues against the other. The arguments we detail bellow can teach us about the limitations of both approaches, but once we trace possible motivations for NF, several similarities become apparent.

Quine is not pleased \cite{QuineSetTheoryAndItsLogic} with Poincar\'e's and Russell's view of the paradoxes, and with the cumulative hierarchy approach. According to the conception advocated by Poincar\'e, Russell, and others, the paradoxes are a result of ``vicious circles'' and ``impredicative definitions'' i.e., definitions of sets with variables that may range over the defined set. Quine explains why there should be no problem with the kind of specifications Poincar\'e warns about and calls impredicative: ``A circular argument seduces its victim into granting a thesis, unawares, as a premiss to its own demonstration. A circular definition smuggles the definiendum into the definiens ... But impredicative specification of classes is neither of these things.'' Thus ``Impredicative specification is not visibly more vicious than singling out an individual as the most typical Yale man on the basis of averages of Yale scores including his own.'' Quine stresses that this interpretation seems to imply an unintended temporal view of sets: ``For we are not to view classes literally as created by being specified---hence as ... increasing in number with the passage of time. Poincar\'e proposed no temporal implementation of class theory. The doctrine of classes is rather that they are there from the start.'' He concludes that ``... the ban urged by Russell and by Poincar\'e is not  to be hailed as the exposure of some hidden ... fallacy that underlay the paradoxes. Rather ... to thin the universe of classes down to the point of consistency.'' 

It seems that Quine's objection to a hierarchical conception of sets stems from his endorsement of definitions that quantify over all sets. According to Quine only definitions of mathematical concepts and deductions should be organized in an expository manner - ``first thing's first.'' A reference to mathematical objects in our domain is something else. Boolos rather claims, following Dan Leary that: \cite{boolos:1989} ``the metaphor of formation of sets at stages may arise from a certain {\em{narrative}} convention or principle of good exposition ...'' Yet Quine seems to think that an expository approach is incompatible with the way we perceive and refer to mathematical objects.

On the other hand, Boolos voices \cite{BoolosIterativeConceptionOfSet} the expected criticism of NF for seeming artificial, and lacking a simple motivation independent of the paradoxes. Referring to Quine's systems NF and ML, Boolos says that: ``Other theories of sets, incompatible with ZF, have been proposed. These theories appear to lack a motivation that is independent of the paradoxes in the following sense: they are not, as Russell has written, `such as even the cleverest logician would have thought of if he has not known of the contradictions.' \cite{Russell1959-RUSMPD-2}'' Boolos points out that Quine's NF lacks a clear motivation, and appears to be artificial. A question mark hovers over NF, because the reasons for accepting only stratified formulas are not clear, which makes NF look like a cleverly ad-hoc way to grab as much as possible while eliminating what we know to be paradoxical just because we are aware of the problem.

Quine does not detail his philosophical considerations for NF in \cite{QuineNewFoundationsForMathematicalLogic}. Later \cite{journals/jsyml/Quine38a} he provides a way to reach his theory through the formalism of type theory. Stratified formulas can be seen as formulas of the theory of types for which the types of the variables have been removed, which is a form of what is known as ``typical ambiguity,'' \cite{Specker1966116} which Quine mentions in that article. This results in a theory with a single type in which what is true within each type of Russell's theory of types (such as having a ``universal set'' relative to each type), becomes true for the single type of NF. 

Quine states his motivation as a way to avoid set antinomies. After examining problematic properties of sets, he finds fault in unstratified formulas and states within this process that: ``We are trying to purify our language of idioms which might deceive us into contradicting ourselves.'' Then he argues that since stratified formulas are formulas of type theory without the types, this should help keep NF free of contradictions, for: ``if the theory of types is adequate at all as a safeguard against contradictions, it must be adequate in its formal aspect alone.'' One may wonder how is that different from an attempt to ``thin the universe of classes down to the point of consistency.''

If Quine performs a mere clever formal trick in that article, Boolos would be right about Quine's approach. Following Ernst Specker \cite{Specker1966116}, it turns out that typical ambiguity has more to it than meets the eye at first glance. It is applicable to other fields of mathematics, and it is relevant to the consistency of NF \cite{SpeckerForthcoming-SPED}. But a motivation for NF which is not formal, and free from complicated mathematical notions has to be sought elsewhere. Holmes \cite{HolmesAdHoc} wrote an article in defense of NF as a framework for set theory. However, the suggested philosophical basis Holmes provides is hardly intuitive, involving notions such as abstract data types when the perception is roughly stated, and ``functions,'' ``machines,'' ``addresses,'' and ``programs'' when the idea is made more explicit.

A more intuitive and graphic suggestion by Sharlow \cite{Sharlow}, which is an extension of the iterative conception, is supposed to support NF. Sharlow begins with the usual iterative conception, but loosens the demand that all collected sets are formed at earlier stages. Sets (and their singletons that appear in the next stage, etcetera), can now be assumed to exist as part of their formation. Thus, for example, a set $S = \{\emptyset, \{S\}\}$ can be formed in the level immediately above the level in which $\emptyset$ is formed, because if we assume $S$ can be created there, then $\{S\}$ would be formed in the next level, and the theory allows collecting sets inspired by $S$ that appear in later stages.

Sharlow explains this kind of set formation by using ``placeholders'' or by thinking all sets exist at once: ``... we can think of the act of placing  $\{S\}$ in $S$ as an act of dropping a {\em{placeholder}} into the set $S$ and making a mental note to replace the placeholder with $\{S\}$ after $\{S\}$ comes into existence. Alternatively, if we want to think of all sets as existing at once, then we can think of ourselves as adding to $S$ an already  existing element which we know must exist once we have created $S$.'' This leaves Sharlow with the need to discuss the difference between self-defeating and non self-defeating properties in his set conception, and show why paradoxical sets are not formed this way. We will encounter the placeholder approach as part of the notion of ``concurrent  aggregation'' we devise for a new set conception in section \ref{sec:complete}.

In summary, the idea of ``typical ambiguity'' shows that NF and the theory of types are not as distant as it might seem at first glance. Also, following Sharlow, a simple motivation for Quine's NF based on the iterative conception is possible. The philosophical rivals that adhere to Russell's type theory and Boolos' stage theory on the one hand, and Quines' NF on the other hand, seem to be closer to each other than they might acknowledge. The layered notions Quine rejected are reflected within the formulas admissible in the comprehension axiom schema of NF. Quine's criticism can now be turned against his own theory. Just like sets, all set properties are supposed to be indifferent to imposed order. How can we justify the rejection of some properties by forcing the variables to be ordered according to their appearance in subformulas? Instead of assuming that elements are prior to their sets, in NF variables that represent elements are prior to variables that represent sets. It looks as though we simply shifted the ordering from sets to variables.

\section{Complete Totalities}
\label{sec:complete}

The question remains if we can transcend this persistent problem of having to accept notions of precedence induced by the membership relation, limitations of size, or any other restriction motivated by the need to maintain consistency. The notion of complete totality we go on to describe is presented as a possible solution to this dilemma.

\subsection{Beyond Rank Theories and Stratified Formulas}

Our objective is to establish a set conception free from the limitations that characterize the conceptions surveyed thus far. As we have seen, both the hierarchical conceptions, and Quine's NF rely on a notion of precedence induced by the membership relation, either between the sets as mathematical objects, or syntactically between variables in subformulas. Otherwise some approaches preserve more basic mathematical intuitions than others. Quine demands and delivers a universe of sets with a single type, and without even a metaphorical notion of precedence of some sets over others. All sets are there ``at once,'' and comprehension quantifies over all sets. Quine's NF supports NWFS including a universal set, though it does not allow the totalities of non-paradoxical unstratified properties such as $x\in x$. The hierarchical approach must include a limitation of size, and reject NWFS. The advantage of cumulative hierarchies is in the assurance we feel about ``aggregating previous sets into new sets or into classes,'' which is not expected to cause contradictions. The ZF approach which can be interpreted as supporting the ``creation'' of sets by both  aggregations and disaggregations, can allow NWFS using anti-foundation axioms instead of regularity, but this still has a ``limitation of size'' limitation. 

What we would like to have is a universe of sets which are there ``at once,'' without any order on sets or variables within formulas. This implies that we are able to quantify over all sets when we define sets, and that there is no ontological dependence between them, so the playground is ``leveled.'' The universe is a universe of sets which follow the notion of a totality as expressed by one of Cantor's definitions for sets that Boolos quotes \cite{BoolosIterativeConceptionOfSet}: ``... a totality of definite elements that can be combined into a whole by a law.'' For the sake of our thought experiment we imagine a procedure of forming the totality of a certain property by first identifying the sets that comply with the property and then aggregating them ``into a whole.'' Of course, none of this carries any physical or temporal notions, similar to the metaphorical speech of ``forming'' sets from ``former'' sets in the iterative conception which does not carry such notions either. Note, however, that there is an inherent conflict between a ``Platonic'' universe of sets and the seemingly safe action of aggregation which entails a conceptual difference between the set as a single unified entity and its elements. Our intention is to take both demands to their extreme, and consider the result, as odd as it might appear at first glance. Insisting on an application of aggregation in a way that avoids even a metaphorical notion of a ``before'' and ``after'' in the way we perceive  aggregation (such as the one expressed by Boolos' ``{\em{Now}} we have a set'') will provide an unexpected way to handle the dilemma.


In order to consider the implications of the conflict mentioned above, it is natural to discuss Russell's property $x \notin x$. The question is of course: does the totality belong to itself or does it not belong to itself? The contradiction stemming from the definition of the totality is many times described as a logical derivation that renders both options symmetric. A careful look shows that things are less symmetric than they seem. If we follow the thought experiment that first identifies all sets such that $x \notin x$ and then forms the totality by aggregating them, we may initially expect that the resulting set does not belong to itself. This may simply arise from the expectation that a collection should be different from each of its elements. It may also arise from the nature of the sets in question, which  do not belong to themselves. Even if we do not expect them to be well-founded, we expect their property to propagate. A confirmation that this is an impossibility comes from the observation that otherwise one of the sets we intended to aggregate must have belonged to itself, thus originally ``misidentified'' and somehow included in the objects to be aggregated.

The asymmetry becomes evident when we notice that the only way to reach the conclusion that the totality belongs to itself is by using a logical derivation based on the property of the totality. The usual argument is: If we assume that the resulting set does not belong to itself, it must be a member of the totality. That is, we apply the creation rule of the collection to itself. This actually has some merit because it delivers on our promise to avoid any notion of ``before'' and ``after'': if the resulting set does not belong to itself, why should it be exempt from its own implications? This would in fact hold also for the singleton of the resulting set and many other sets that can be theorized using it.

To put it in a way that conforms with the ``placeholder'' approach mentioned above, if the resulting set does not belong to itself, we have again ``misidentified'' a set, this time by forgetting to identify the resulting set as one of the original elements that comply with the property, thus we are forced to assume that we originally ``had something there'' and got an instance of what we are collecting. According to this view, at the instant of collecting, the collected identifies with the collection, so the end result will have itself as an element. We call this ``concurrent aggregation.''

There is another way to describe the same effect. We initially expect to get something  which is different than the the sets we try to aggregate, all of which do not belong to themselves. Hmmmm, it looks like we ``forgot'' one. Not discouraging, since we assume there is no ``before'' and ``after'' in our universe, we ``force it back'' into the initial group of sets we try to collect, and try to aggregate ``again,'' only to get another one, and another one, ad infinitum. Each time we aggregate we get another, bigger, set, and the descending $\in$-sequence of these sets gets one step longer. At the ``limit'' of this process we get an infinitely descending sequence of ... indistinguishable sets! Assuming a strong notion of equating NWFS, the result is a set belonging to itself. According to similar reasoning, the singleton of this set, and many other sets that contain it in some way, must be members of the resulting set.

What did we get? Well, it is obvious what we did not get. This cannot be the totality of $x \notin x$ because it includes sets that do not comply with the property $x \notin x$. On the other hand, the creation of this set does not seem arbitrary, and unlike the previous situation in which we were unable to respond to the question mark induced by $x \notin x$, we now have something we might want to investigate further. In order to distinguish the two collections we will call the usual totality an {\em{ideal totality}} of a property, which is known not to be a set in some cases, and the new collection the {\em{complete totality}} of a property, which we intend to assert as an axiom to always be a set. Several immediate questions follow: Is there a way to characterize these seemingly foreign intruders into the totality? What could possibly be their meaning? Does the emerging universe of sets seem coherent, what does it contain, and can it be useful?

\subsection {Forming Complete Totalities}

Let us try to clarify the idea of a complete totality. The general idea of a complete totality is that instead of just acknowledging all sets that comply with the predicate, we also seek the implications of their collection. If a hypothetical object, not necessarily a set, created from the ideal totality and other sets (that is, if $I$ is the ideal totality we will consider structures such as $I, \{I\}, \{I,\emptyset\}, \{\{I\}\}$ etcetera) complies with the predicate, this would indicate the existence of an additional element in the complete totality. However, that additional element is not the above hypothetical object, but the set that results from the hypothetical object for which all the occurrences of the ideal totality are replaced by the complete totality (that is, if $C$ is the complete totality we respectively accept $C, \{C\}, \{C,\emptyset\}, \{\{C\}\}$ as members of $C$). By ``replace'' we mean something like using the occurrences of the ideal totality as ``placeholders,'' the way placeholders were described in Sharlow's broadening of the iterative conception we encountered earlier. Unlike Sharlow though, we do not keep stages, and we allow the complete totality of all predicates. The intention is to adopt as an axiom the assertion that all complete totalities are sets, unlike some ideal totalities.

This concept sounds somewhat complicated, but once the idea of replacement using placeholders is made clear, all the definitions become simple. We will use an informal discussion and ad-hoc notation, including the vague notion of ``hypothetical object,'' because the intention is only to apply them to important predicates and make some immediate observations. An informal definition of replacement would be as follows:


\begin{definition}
The {\em{replacement}} of $x$ with $y$ within $S$, denoted by $S_{x\rightarrow y}$ is:
    \begin{enumerate}

        \item If $S=x$, then let $S_{S\rightarrow y} = y$.
        
        \item Otherwise $S_{x\rightarrow y} = \{ z_{x\rightarrow y} | z\in S\}$.

    \end{enumerate}

\end{definition}

For example, this means, assuming the usual meaning of 2 and $\omega$, that $2_{\emptyset\rightarrow \omega} = \{\emptyset, \{\emptyset\} \}_{\emptyset\rightarrow \omega}  = \{\omega, \{\omega\} \}$. Even though we know that many more instances of $\emptyset$ appear within $\omega$, those will remain as are and not be replaced, that is, the replacement is never propagated. Of course, this kind of ``definition'' only works for well-founded sets, since according to this there is no way to determine $\Omega_{\emptyset\rightarrow\omega}$ where $\Omega$ is the singleton of itself. We will take this definition to mean that if along an infinitely descending $\in$-sequence we encounter the element to be replaced, it is replaced the first time it is noticed down the sequence. Otherwise the element of be replaced does not belong to the transitive closure of the object replacement operates on, so replacement becomes latent, and the end result is identical with that object. Thus for example $\Omega_{\emptyset\rightarrow\omega} = \Omega$, and $\Omega_{\Omega\rightarrow\omega} = \omega$.

Moving on to the two notions of totality, we have the following:

\begin{definition}
The {\em{ideal totality}} of a predicate $P$, denoted by $[P]$ is the collection of all sets that satisfy $P(x)$.
\end{definition}

\begin{definition}
The elements of a {\em{complete totality}} of a predicate $P$, denoted by $(P)$ are either:

    \begin{enumerate}

        \item All sets $x$ such that $P(x)$ holds, are elements of $(P)$.
        
        \item For any non-set $x$ such that $P(x)$ holds: if $x_{[P]\rightarrow (P)}$ is a set, it is an element of $(P)$.

    \end{enumerate}

\end{definition}

One might say that the definition of a complete totality is a bit artificial, because it has two branches. However, both branches go by the same reasoning - if something satisfies $P$, ``before'' or ``after'' we aggregate, it indicates an element that should be included. The complete totality simply treats both possibilities as equivalent. 


There are two problems with this ad-hoc definition of complete totalities. The first is that it does not include a way to determine how many non-sets that include $[P]$ in their transitive closure we have, that is, it keeps the criteria for ``hypothetical objects'' that can be built using $[P]$ implicit. The second problem is that this definition is circular. We will see that this circularity will make the definition depend on aspects of our set theory that may be decided in different ways, such as the notion of set equality. However, it turns out that many things can be decided about complete totalities even with such a simple exposition.

We are ready to state the axiom of complete totalities:

\begin{axiom}

For any property $P$, $(P)$ is a set.

\end{axiom}

This might be seen as an extension of Quine's NF, once there is also an axiom of extensionality. However, the analysis of some simple properties we shall soon discuss may lead to two conclusions that would make such a system different from NF in aspects other than allowing unstratified formulas. The first one is that it might be a good idea to use a stronger version of extensionality, which is more suitable for work with non-well-founded sets. The second one is that unlike NF which allows formulas that may contain more than one free variable (though, as usual with various comprehension axioms, the variable that indicates the set formed by the axiom is rejected in NF), that the properties we want to allow will contain only a single free variable.


What kind of a universe of sets do we get? An application of the concept to several important properties will deliver some initial answers. Obviously $(x\neq x) = [x\neq x] = \emptyset$. This must be because any route of acceptance of elements into the complete totality requires that something satisfies the property in question.

A quick reflection shows that it is easy to create a property (we consider a property to be at least everything that can be expressed by first order logic formulas in the language of set theory) that allows any specific well-founded\footnote{Aczel uses the term hereditarily finite in a way that includes non-well-founded sets.} hereditarily finite set. The ideal totality in this case and anything built using it are not expected to satisfy the property. 

An infinite set, in fact $\omega$ itself, is created using the property for finite ordinals. Since the ideal totality in this case is infinite, and all finite ordinals are well-founded hereditarily finite sets, $\omega$ will also be the complete totality. Many more ordinals and sets required for ``everyday mathematics'' with well-founded elements of bounded rank are easily obtained in a similar way.

The universal set can be obtained using $x=x$. Obviously $(x=x) = [x=x]$. Note that we already know that this universal set must be infinite, and contain infinite sets. It must also contain itself, the singleton of itself, and many other variants which keep it in their transitive closure. What we get seems to depend on the notion of a ``hypothetical object'' which remains obscure, but even if we allow only all structures that can be built using a finite number of aggregations from $[x=x]$ and all other sets, the number of combinations is overwhelming. This makes our universal set seem like a universal set worthy of its name. That is, it is not simply a set that satisfies the formula $\forall x x\in y $, which is satisfiable even in a model with a domain of size one.

The property $Ord$ of ordinals is interesting both because it is possible to characterize a complete totality which is different than the ideal totality even when ``possible structure'' remains obscure, and also in terms of the implications of different kinds of extensionality. Since we know that in this case the ideal totality is an ordinal, apart from all the ordinals the complete totality will include itself as an element. Can there be any other elements? Since the successor of the ideal totality is also an ordinal, replacement will yield a set which includes all the ordinals and the complete totality. This seems to be exactly the same as what we already have, so we may expect that all further ``successors'' and ``limit ordinals'' (all of which are not sets), must yield an object similar to the complete totality, which is already included within the complete totality. It looks as though we can conclude that the complete totality will keep only the ordinals and itself.

However, if we do the above symbolically we see that identifying the objects created by replacement with the complete totality depends on our notion of set equality. Let us try this exercise carefully. For simplicity let $\alpha = [Ord]$ and $\beta = (Ord)$. Since $Ord(\alpha)$ and $\alpha_{\alpha\rightarrow \beta} = \beta$, and $\beta$ is a set, we have $\beta\in\beta$ thus $\alpha \cup \{\beta\} \subseteq \beta$. We expect that for all ordinals $\gamma$, $(\alpha + \gamma)_{\alpha\rightarrow \beta} = \beta$, and we expect that this simply propagates, so the standard way would be to use transfinite induction. We assume for all $\delta<\gamma$ that $(\alpha + \delta)_{\alpha\rightarrow \beta} = \beta$. Thus, $(\alpha + \gamma)_{\alpha\rightarrow \beta} = (\alpha \cup (\cup_{\delta<\gamma}{\{\alpha+\delta\}}))_{\alpha\rightarrow \beta} = \alpha \cup (\cup_{\delta<\gamma}{\{\beta\}}) = \alpha \cup\{\beta\}$. The only step missing is $\alpha \cup\{\beta\} = \beta$, but we do not have it. What we have is $\alpha \cup \{\beta\} \subseteq \beta$. Our analysis shows that once the identity occurs somewhere it would propagate. Otherwise we get within $\beta$ a set that our identity relation does not equate with $\beta$, composed of itself and all ordinals. Or we may get that same situation a finite or an infinite times down an $\in$-sequence, depending the identity relation, and the specific formalization and model. However, fixing a strong sense of equality between sets (for example by using the anti-foundation axiom with the strongest sense of equality described by Aczel) to force $(Ord)$ to be a set containing all ordinals and itself. 

Interestingly, complete totalities make a clear distinction between collecting all finite ordinals, and collecting all ordinals, based on the fact that the collection of finite ordinals is infinite, but the collection of all ordinals is an ordinal. It could be said that complete totalities have their own natural ``size considerations mechanism'' which acts upon principle and not upon arbitrary decisions.

Some things are expected not to be the way they are for mainstream systems such as ZF. Since we have a universal set, just like in NF, we expect that either some sets do not have a power set, or that it impossible to prove Cantor's theorem about power sets always having greater cardinalities. We also know that some aspects of the system we describe must be more peculiar than NF. For example, in NF each set has a complement, because the negation of a stratified formula is stratified. However we know that $(x\in x)$ and $(x\notin x)$ are not disjoint, since $(x\notin x)$ belongs to both sets.

What can we expect to have in terms of power sets, union sets, and other methods of ``creating'' sets from other sets? An easy exercise shows that  there is a limitation on the properties we can use. Unlike NF, the formulas for the properties we use for the elements of sets cannot include free variables other than the one we use for the elements themselves. NF allows parameters inside the properties through free variables that do not include the variable for the defined set. This is useful for getting power sets, and performing other set operations. In our case, however, using the notations we used in the discussion about ordinals, $\beta \smallsetminus \{\beta\}$ would yield a set $\alpha$. In formula notation this would happen if we allow the formula $x\in y \wedge x\notin x$, in which $x$ describes the elements of the totality, and $y$ is free. If this is allowed the axiom of complete totalities would mean that there is a complete totality (thus a set) for the property, whenever $y$ is assigned a given set. Since $\beta$ is a set, the property for $\beta$ assigned to $y$ is about all $x\in\beta$, such that $x\notin x$, and those are the ordinals. This time the ideal totality of the property, namely $\alpha$, and all hypothetical structures containing it, do not comply with the property because they are not in $\beta$.

What is the effect of this limitation? Since our original observations were made for properties of sets and not for properties of sets that have parameters - because the basic premise was that we should be able to aggregate any group of well defined sets - there is no dent in the original conception. On a practical level it means that there is no direct and usual way to assert the existence of sets, given the existence of other sets. NF always allows a power set this way, and blocks Cantor's paradox by not allowing the diagonal set needed for the theorem, and by other conditions that arise when similar sets with a diagonal set are considered. In our case we simply do not have a power set for all sets, but power sets for important sets such as $\omega$ are easily created using direct properties. In this case we simply try to create the totality of all sets that have only finite ordinals as members, which are equivalent of subsets of $\omega$. Luckily properties of well-founded sets with a limited rank lend themselves easily to our procedures. The ideal totalities and the complete totalities become identical in such cases. It seems that even with this limitation we are not expected to have trouble defining the sets required for ``usual mathematics.''

The picture is getting clearer, and at the same time many questions and possibilities arise. It would be natural to explore other predicates, to offer a mathematical formalization that handles circularities in a rigorous way and would make all notions precise, to prove that the formalization is consistent relative to ZF, to seek the implications on large cardinals, and on theorems that cannot be decided using other conceptions of sets, etcetera. Since the focus of this article is on presenting a new set conception and on its philosophical advantages with respect to other set conceptions, we leave the above issues for further research. 

\subsection{Amiable Circles and the Virgin Paradox}
 
We finish with a few philosophical observations on the notion of complete totality, which include a consideration of the legitimacy and meaning of elements that do not comply with a property yet reside in the property's complete totality, and the implications of that meaning on the property. Also, a short note about the philosophical implications of considering a complete totality to be the ``official totality'' of a property.

As for the meaning of the additional elements, note that elements of a complete totality that occur as an act of (non latent) replacement, either ones that eventually comply or ones that do not comply with the property, are always circular, and all of them occur as a result of contemplating the implications of having a the usual ideal totality. The meaning is therefor pointing out circularities within the property itself. For example, the fact that the complete totality of the ordinals keeps itself as an element, indicates that the notion of ``all ordinals'' is in itself an ordinal. This valuable information about the property becomes embedded into its complete totality. Reflecting upon Russell's ``viscous circle principle'' which states that ``no totality can contain members defined in terms of itself.'' \cite{Rus:08} we note that the notion of a complete totality not only rejects this principle, but actively does define members in terms of the totality in the more problematic cases, which makes the inherent circularities of the property visible. Instead of adopting an alarming title and fighting these circularities, one can accept their natural emergence and investigate them. If titles are of any interest, those can be ``amiable circles'' just as well.

The legitimacy of keeping elements that do not comply with a property in its totality was based on the fact that in order for anything to be included in a complete totality, the property must apply to some object, even if it is only hypothetically based on the contemplation of the ideal totality. This form of set creation was said to  explore the implications of aggregating all sets that have a certain property, and act upon them. It could be argued that this is just a form of mathematical or conceptual wizardry. How can we see as natural the situation in which an operator must create a situation that negates the meaning of its operands?

A story titled ``the virgin paradox'' might make this clear: Randolph Higgins, better known to his friends as ``randy Randolph'' for his long history of boasting his sexual performance, sits in a bar with his friends and claims he had sex with a virgin earlier that week. None of the participants sees this statement as strange, and they urge him to deliver details about that incident. This is not a result of them knowing Mr. Higgins. None of us would find a statement about having sexual intercourse with a virgin (which literally means a woman that never had sex with anyone) strange, because even though sex and virgins are somewhat contradictory concepts, the meaning of the statement is perfectly clear. These concepts are mutually exclusive in the semantic sense, but not mutually exclusive as an actual occurrence. To have sex with a virgin must be an instant occurrence that at the same time changes the status of the woman. There is no way to have sex with a virgin if she already had a man, for she would not be a virgin. There is no way to have sex with a virgin, and have her change her status at some later time. There is even no exact way to define the instant of the transformation. This kind of an exact point which is also a blur could almost be a good proof for the impossibility of sexual reproduction. Yet we know that sex and virgins mix on a daily basis.

The situation of aggregations, which we mentally perceive to be some kind of operator, on the one hand, and of a universe of mathematical objects that we prefer to see as static and unchanging on the other hand, is somewhat similar to the virgin paradox. This is probably the deeper reason behind set antinomies.

Finally, a note about collectives and meaning. Complete totalities have many merits. They exist for any property, they make hidden circularities evident, and they are a result of aggregation that takes into account its own implications. Can anyone claim that complete totalities are ``really the'' totalities of all properties, even when they may contains members that do not comply with the property attached to the same totality? This appears to be a strange philosophical proposition, which is not required from anyone interested in the purely mathematical aspects of the axiom of complete totalities. An adoption of complete totalities as ``true'' totalities would imply that plurals we address as a single unified entity must have a meaning that transcends the meaning of the predicate that instigated them. This would be immediately rejected by anyone who perceives sets as nothing more than a combination of their elements, and even more so by anyone claiming that a set is the same as its elements\footnote{For a support of the latter claim see for example \cite{boolos:1989}.}. However, it is possible to argue that not only collectives may have an independent meaning, but that any collective that may be of any interest must have one. This argumentation is beyond the scope of this article.

In summery, in this article we have presented a new notion of sets and of totalities, complete with its basic motivation, and some of its immediate mathematical implications. As a philosophical position it offers a new outlook on the reasons for set theoretic paradoxes. It also delivers a universe of sets close in nature to the naive conception of set. As a mathematical suggestion, it might be able to offer some new insights about large cardinals, and about important set theoretic postulates known to be undecidable in ZF. This is left for further research.

\bibliography{CompleteTotalities} 

\vspace{ 8mm}
Copyright \textcopyright \hspace{1mm} 2011 Rafi Shalom.

\end{document}